# Lower bounds for posterior rates with Gaussian process priors

**Ismaël Castillo**

*Department of Mathematics,*
*Vrije Universiteit Amsterdam,*
*De Boelelaan 1081a, 1081HV Amsterdam, Nederland.*
*e-mail:* i.castillo@few.vu.nl

**Abstract:** Upper bounds for rates of convergence of posterior distributions associated to Gaussian process priors are obtained by van der Vaart and van Zanten in [14] and expressed in terms of a concentration function involving the Reproducing Kernel Hilbert Space of the Gaussian prior. Here lower-bound counterparts are obtained. As a corollary, we obtain the precise rate of convergence of posteriors for Gaussian priors in various settings. Additionally, we extend the upper-bound results of [14] about Riemann-Liouville priors to a continuous family of parameters.

**AMS 2000 subject classifications:** Primary 62G05, 62G20.
**Keywords and phrases:** Bayesian nonparametrics, Gaussian process priors, Lower bounds.

## 1. Introduction

In the Bayesian non-parametrics literature, several general results about posterior consistency (see e.g. [1]) and posterior rates of convergence (see for instance [5; 13]) are now available. Roughly, the rate of convergence of the posterior is generally thought of as an $\varepsilon_n$ as small as possible such that the posterior probability of the ball centered at the true $f_0$ and of radius $\varepsilon_n$ still tends to 1 in probability. In this context a natural question is, starting from a fixed prior, what is the actual rate of convergence of the posterior ? The tools proposed in the cited articles often allow to get an upper bound for this posterior rate.

Also, from the practical point of view, non-parametric type priors are now commonly used in applications, as an example the book [12] presents applications of Gaussian priors in machine learning. In non-parametric situations many priors will not lead to optimal rates; in some cases the corresponding posterior will still converge at some reasonable rate towards the true parameter or function; in other cases the convergence might be extremely slow or consistency might even fail. Determining the precise rate of convergence of the posterior can then help in choosing the type of prior adapted to the practical situation or in adjusting the prior parameters.

Given a class of functions, upper bounds for the rate are clearly optimal if they coincide with the minimax rate of convergence over the class. In the case where the upper bound is slower than the optimal rate, one would like to establish a bound from below for the rate. If both upper and lower bounds match, say up to some constant or logarithmic factor, then the exact rate of convergence of the posterior remains determined. In this paper, the issue of obtaining a lower bound for the posterior rate is considered in the case of Gaussian priors. Though the focus will be mainly on the class of Gaussian priors, the methodology we introduce can be used in some cases to obtain lower bounds for other priors as well. In particular, we shall derive a lower bound result for a non-Gaussian prior in a specific example.

The organization of the paper is as follows. In the next section, we enounce our main result on lower bounds in a general framework and give its proof. This result is applied in Section 3 to obtain lower bounds in two nonparametric models: the Gaussian white noise model and the problem of density estimation, with respectively Gaussian series priors and Riemann-Liouville priors. For the latter prior, upper bounds are also obtained which extend previous results of [14]. Technical results are gathered in Section 4. Concluding remarks are given in Section 5.

Let us introduce some notation. For any real numbers $a, b$, we denote by $a \wedge b$ their minimum and by $a \vee b$ their maximum. We define Hellinger's distance $h(f, g)$ between two probability densities $f$ and $g$ by the $L^2$-distance between the root densities $\sqrt{f}$ and $\sqrt{g}$. Let $K(f, g) = \int f \log(f/g) d\mu$ stand for the Kullback-Leibler divergence between the two non-negative densities $f$ and $g$ relative to a measure $\mu$. Furthermore, we define the additional discrepancy measure $V_2(f, g) = \int f |\log(f/g) - K(f, g)|^2 d\mu$. Let $L^2[0, 1]$ be the space of square integrable functions on the interval $[0, 1]$, equipped with the $L^2$-norm $\|f\|_2 = (\int_0^1 f^2 d\mu)^{1/2}$. Let $\mathcal{C}^0[0, 1]$ denote the space of continuous functions on $[0, 1]$ equipped with the supremum norm $\|\cdot\|_\infty$. Let $\mathcal{C}^\beta[0, 1]$ denote the Hölder space of order $\beta$ of all continuous functions $f$ that have $\underline{\beta}$ continuous derivatives, for $\underline{\beta}$ the largest integer strictly smaller than $\beta$, with the $\underline{\beta}$th derivative $f^{(\underline{\beta})}$ being Lipshitz-continuous of order $\beta - \underline{\beta}$. This means that there exists a positive constant $C$, which might depend on the function $f$, such that

$$|f^{(\underline{\beta})}(y) - f^{(\underline{\beta})}(x)| \leq C|y - x|^{\beta - \underline{\beta}} \quad \forall \, x, y \in [0, 1].$$

## 2. Lower bound result

Let $(\mathcal{X}^{(n)}, \mathcal{A}^{(n)}, P_f^{(n)}; f \in \mathcal{F})$ be a sequence of statistical experiments with observations $X^{(n)}$, where the parameter set $\mathcal{F}$ is a subset of a Banach space $\mathbb{B}$ (for instance $L^2[0, 1]$ or $C^0[0, 1]$) and $n$ is an indexing parameter, usually the sample size. We put a prior distribution $\Pi$ on $f$. In this paper we consider the case where the prior is the law of a Gaussian process taking almost surely its values in $\mathbb{B}$ (see below). We are interested in properties of the posterior distribution $\Pi(\cdot|X^{(n)})$ under $P_{f_0}^{(n)}$, where $f_0$ is the "true" function. We denote by $\mathbf{E}_0$ the expectation under the latter distribution. For any sequence $\varepsilon > 0$ let us define a Kullback-Leibler neighborhood of $f_0$ as

$$B_{KL}(f_0, \varepsilon) = \{f: \ K(P_{f_0}^{(n)}, P_f^{(n)}) \leq n\varepsilon^2, \ V_{2,0}(P_{f_0}^{(n)}, P_f^{(n)}) \leq n\varepsilon^2\}.$$

In this work Gaussian processes $Z$ are supposed to be centered and tight measurable random maps in the Banach space $(\mathbb{B}, \|\cdot\|)$. We refer to [15] for an overview of basic properties of these objects. Let $\mathbb{H}$ be the Reproducing Kernel Hilbert Space (RKHS) of the covariance kernel of the process. We will generally assume that $f_0$ belongs to the support of the prior in $\mathbb{B}$, which for Gaussian process priors is nothing but the closure of $\mathbb{H}$ in $\mathbb{B}$, see for instance [15], Lemma 5.1.

First let us review the key upper-bound results obtained in [14], where the authors show that for Gaussian priors an upper-bound for the concentration rate of the posterior can often be obtained in a simple way from the so-called *concentration function* of the Gaussian process. This quantity is defined as follows. For any $\varepsilon > 0$, let

$$\varphi_{f_0}(\varepsilon) = \inf_{h \in \mathbb{H}: \|h - f_0\| < \varepsilon} \|h\|_{\mathbb{H}}^2 - \log \mathbf{P}(\|Z\| < \varepsilon) \tag{1}$$

Assume that the norm $\|\cdot\|$ on $\mathbb{B}$ is comparable to a metric $d$ appropriate to the statistical problem (often, $d$ is a distance for which certain tests exists, which allows to apply the theory presented in [5]; for instance, in i.i.d. settings, one might choose Hellinger's distance). Here "comparable" means that the ball $\{f \in \mathcal{F}, \|f - f_0\| \leq \varepsilon_n\}$ should be included in the ball for $d$ around $f_0$ of radius $c\varepsilon_n$ and also in the Kullback-Leibler neighborhood $B_{KL}(f_0, c\varepsilon_n)$ defined above, for some $c > 0$. The authors in [14] prove that if $\varepsilon_n \to 0$ satisfies

$$\varphi_{f_0}(\varepsilon_n) \leq n\varepsilon_n^2, \tag{2}$$

then the posterior contracts at the rate $\varepsilon_n$ for the distance $d$, in that for large enough $M > 0$, $\mathbf{E}_0 \Pi(f: d(f, f_0) \leq M\varepsilon_n \mid X^{(n)}) \to 1$ as $n \to \infty$.

These results mean that for Gaussian priors an upper-bound on the rate of the posterior is obtained as soon as the next two quantities are controlled

$$\varphi_{f_0}^A(\varepsilon) = \inf_{h \in \mathbb{H}: \|h - f_0\| < \varepsilon} \|h\|_{\mathbb{H}}^2, \qquad \varphi^B(\varepsilon) = -\log \mathbf{P}(\|Z\| < \varepsilon). \tag{3}$$

The first term measures how well elements in the RKHS $\mathbb{H}$ of the Gaussian process can approximate the true function. Note in particular that if $f_0$ happens to be in $\mathbb{H}$, this term simply remains bounded. The second term, which does not depend on $f_0$, is the so-called *small ball probability* of the Gaussian process. Small ball probabilities have been studied in many papers in the probability literature and precise equivalents as $\varepsilon \to 0$ of $\varphi^B(\varepsilon)$ are available for many classes of Gaussian processes, see for instance [11]. Yet at first sight it is not obvious to see why the concentration function $\varphi_{f_0}$ should appear in the study of posterior rates. Lemma 2 below answers, at least partially, this question.

Let us conclude the overview of upper-bound results with an example. In a context of density estimation, if one chooses Brownian motion as prior on continuous functions, the rate $\varepsilon_n$ depends on the Hölder regularity $\beta$ of the true $f_0$ as follows. If $\beta \geq 1/2$, then $\varepsilon_n$ can be chosen equal to $n^{-1/4}$, whereas if $\beta < 1/2$ the rate $\varepsilon_n$ must be in $n^{-\beta/2}$ to satisfy (2), see Section 4.1 in [14] or Theorem 3 below. Thus, up to constants, the rate is optimal in the minimax sense if $\beta = 1/2$. However, for all other values of $\beta$, the obtained rate is below the minimax rate which is $n^{-\beta/(2\beta+1)}$. Thus it is natural to ask whether the rate of concentration for Brownian motion is really the one described above or if in fact the posterior contracts faster.

Let us now define a notion of lower bound for a given prior $\Pi$ on $\mathcal{F}$. Let $d$ be a distance on the parameter space. We say that the rate $\zeta_n$ is a lower bound for the concentration rate of the posterior distribution $\Pi(\cdot|X^{(n)})$ in terms of $d$ if, as $n \to +\infty$,

$$\mathbf{E}_0 \Pi(f: d(f, f_0) \leq \zeta_n \mid X^{(n)}) \to 0. \tag{4}$$

This mainly means that $\zeta_n$ is too fast for the posterior measure to capture mass in the ball of radius $\zeta_n$ around $f_0$. The aim is then to prove that such a result holds for $\zeta_n$ as large as possible. In the sequel, we will be able to prove in some examples that the posterior puts asymptotically all its mass inside a *ring* of the type $\{f \in \mathcal{F}, m_n \varepsilon_n \leq d(f, f_0) \leq \varepsilon_n\}$, for $m_n$ either a small enough constant or slowly decreasing (e.g. of logarithmic order), see Section 3. Note also that a lower bound in the sense of Definition (4) will not be an upper-bound for the same distance, so our definition is, in a way, in a strict sense. But it seems to us to be the most natural one, for symmetry reasons with respect to upper-bound definitions, and also in view of the aforementioned 'ring'-behavior. It would also be interesting to be even more precise about the behavior of the posterior: for instance, if asymptotically the posterior sits on a ring for a distance $d$, to see how the mass is distributed inside this ring. However, the presently available techniques, including the ones of this paper, give only results up to constants, so this would probably require introducing new techniques or refining the mentioned ones.

Theorem 1 below establishes a lower bound for the concentration rate of the posterior $\Pi(\cdot|X^{(n)})$ for Gaussian priors in terms of the norm $\|\cdot\|$ of the Banach space. Its proof relies on two basic ideas. The first one is that, roughly, if the prior probability puts very little mass (in some sense) on a certain measurable set, then the posterior probability of this set is also small. The following lemma is Lemma 1 in [6] (see also Lemma 5 in [1]).

**Lemma 1.** *If $\alpha_n \to 0$ and $n\alpha_n^2 \to +\infty$ and if $B_n$ is a measurable set such that*

$$\Pi(B_n)/\Pi(B_{KL}(f_0, \alpha_n)) \leq e^{-2n\alpha_n^2},$$

*then $\mathbf{E}_0 \Pi(B_n \mid X^{(n)}) \to 0$ as $n \to +\infty$.*

The second ingredient is a general result about Gaussian priors which gives control from above and below of non-centered small ball probabilities associated to the process in terms of $\varphi$. For a proof, see for instance [9] or [15], Lemma 5.3.

**Lemma 2.** *Let $Z$ be a Gaussian process in $\mathbb{B}$ with associated RKHS $\mathbb{H}$. Assume that $f_0$ belongs to the support of $Z$ in $\mathbb{B}$. Then for any $\varepsilon > 0$,*

$$\varphi_{f_0}(\varepsilon) \leq -\log \mathbf{P}(\|Z - f_0\| < \varepsilon) \leq \varphi_{f_0}(\varepsilon/2).$$

In view of this result, it seems natural to see the concentration function $\varphi_{f_0}$ appear when studying rates of contraction for Gaussian processes, since the latter function gives a direct control on how much mass the prior puts on neighborhoods of the true function. The following lemma now states some useful properties of $\varphi_{f_0}$. In particular, it implies that this function has an inverse $\varphi_{f_0}^{-1}$.

**Lemma 3.** *Let $Z$ be a non-degenerate centered Gaussian process in $(\mathbb{B}, \|\cdot\|)$. For any $f_0$ in $\mathbb{B}$, the associated concentration function $\varepsilon \to \varphi_{f_0}(\varepsilon)$ is strictly decreasing and convex on $(0, +\infty)$. In particular, it is continuous on $(0, +\infty)$.*

This lemma is proved in Section 4. We can now state our first Theorem.

**Theorem 1.** *Let $Z$ be a Gaussian process with associated distribution $\Pi$ on the space $(\mathbb{B}, \|\cdot\|)$. Let the data $X^{(n)}$ be generated according to $P_{f_0}$ and assume that $f_0$ belongs to the support of $\Pi$ in $\mathbb{B}$. Let $\alpha_n \to 0$ such that $n\alpha_n^2 \to +\infty$ and $\Pi(B_{KL}(f_0, \alpha_n)) \geq \exp(-cn\alpha_n^2)$ for some $c > 0$. Suppose that $\zeta_n \to 0$ is such that $\varphi_{f_0}(\zeta_n) \geq (2 + c)n\alpha_n^2$. Then, as $n \to +\infty$,*

$$\mathbf{E}_0 \Pi(\|f - f_0\| \leq \zeta_n \mid X^{(n)}) \to 0.$$

*Proof.* Due to Lemma 2, it holds $\Pi(\|f - f_0\| \leq \zeta_n) \leq \exp(-\varphi_{f_0}(\zeta_n))$. Combining this with the assumption on the KL-type neighborhood, one gets that

$$\frac{\Pi(\|f - f_0\| \leq \zeta_n)}{\Pi(B_{KL}(f_0, \alpha_n))} \leq \exp(-\varphi_{f_0}(\zeta_n) + cn\alpha_n^2).$$

The assumption on $\varphi_{f_0}(\zeta_n)$ ensures that the last display is further bounded from above by $\exp(-2n\alpha_n^2)$. An application of Lemma 1 with the choice of set $B_n = \{f \in \mathcal{F}, \|f - f_0\| \leq \zeta_n\}$ leads to $\mathbf{E}_0 \Pi(B_n \mid X^{(n)}) \to 0$. □

Before commenting on this result, let us state a direct consequence of it. If the rate $\varepsilon_n$ satisfies (2) and if the norm $\|\cdot\|$ combines correctly with the Kullback-Leibler divergence, so that for some $d > 0$, it holds

$$\Pi(B_{KL}(f_0, d\varepsilon_n)) \geq \Pi(\|f - f_0\| < 2\varepsilon_n),$$

see Section 3 or [14] for some examples, then due to Lemma 2 we obtain that $\Pi(B_{KL}(f_0, d\varepsilon_n)) \geq \exp(-n\varepsilon_n^2)$. Hence according to Theorem 1,

$$\zeta_n = \varphi_{f_0}^{-1}((1 + 2d^2)n\varepsilon_n^2)$$

is a lower bound for the rate of convergence.

Furthermore, if $\varphi_{f_0}$ is "nicely varying" (see below, this depending of course on the particular function $f_0$), then one expects to be able to chose $\zeta_n$ of about the same order as $\varepsilon_n$ (e.g. $\zeta_n = \varepsilon_n/\log n$ or even $\zeta_n = \varepsilon_n/K$ for $K$ large enough). For instance, if $\varphi_{f_0}^{-1}$ is of regular variation in the neighborhood of $+\infty$, then $\zeta_n(f_0)$ is at least $\varepsilon_n/K$, for some $K$ large enough.

Thus we complement the result of [14], where the upper bound part was obtained, by proving a lower bound counterpart. Note also that interestingly, to prove Theorem 1, just the lower bound of Lemma 2 is used. By contrast, note that the main ingredients of the proof of the upper bound in [14] are Borell's inequality and the upper bound of Lemma 2. Note also that the assumptions of Theorem 1 are mainly in terms of the prior, the model coming in only through the Kullback-Leibler neighborhood.

As stated Theorem 1 can be used for Gaussian priors only. However, it illustrates well how the simple Lemma 1 can successfully be applied to obtain lower bounds for the concentration rate. In fact, when dealing with general priors, one can try to apply Lemma 1 directly. This idea enables us to obtain a lower bound result for a non-Gaussian prior (though constructed from Gaussian priors) further in this paper, see Theorem 3. This approach seems to be useful to get lower bounds for general priors in other contexts as well. Further contributions on this question are in preparation and should be available soon.

Another interesting question is how to get more explicit estimates of the rates $\varepsilon_n$ and $\zeta_n$ in terms of the class of functions the true $f_0$ belongs to and of the "regularity" $\alpha$ of the process in some sense (for Brownian motion and Hölder classes we would have $\alpha = 1/2$). In the next section, we address this question in some simple cases.

## 3. Applications

### 3.1. The $L^2$-setting and Gaussian series priors

Let $\{\varepsilon_k\}_{k\geq 1}$ be an orthonormal system in $L^2[0,1]$, being chosen for simplicity equal to the trigonometric basis $\varepsilon_1 = 1$ and for $k \geq 1$, $\varepsilon_{2k}(\cdot) = \cos(2\pi k\cdot)$ and $\varepsilon_{2k+1}(\cdot) = \sin(2\pi k\cdot)$. The Sobolev ball $\mathcal{F}_{\beta,L}$ of order $\beta > 0$ and radius $L > 0$ is

$$\mathcal{F}_{\beta,L} = \{f \in L^2[0,1],\ f = \sum_{k\geq 1} f_k \varepsilon_k \quad \text{and} \quad \sum_{k\geq 1} k^{2\beta} f_k^2 \leq L^2\}.$$

*Gaussian series priors.* Let $\{\alpha_k\}_{k\geq 1}$ be a sequence of independent standard normal random variables and let $\{\sigma_k\}_{k\geq 1}$ be some square-integrable sequence of real numbers. For simplicity let us choose $\sigma_k = k^{-1/2-\alpha}$ for some $\alpha > 0$. Let us define $\Pi$ as the probability distribution generated by

$$X_\alpha(\cdot) = \sum_{k=1}^{+\infty} \sigma_k \alpha_k \varepsilon_k(\cdot). \tag{5}$$

This defines a process with sample paths in $\mathbb{B} = L^2[0,1]$. The RKHS $\mathbb{H}^\alpha$ of $X_\alpha$ in $\mathbb{B}$ is $\mathbb{H}^\alpha = \{\sum_{k\geq 1} h_k \sigma_k \varepsilon_k,\ (h_k)_{k\geq 1} \in l^2\}$, equipped with the norm $\|\sum_{k\geq 1} h_k \sigma_k \varepsilon_k\|_{\mathbb{H}_\alpha}^2 = \sum_{k\geq 1} h_k^2$, see for instance [15], Theorem 4.2. Since the support of the process in $L^2$ is then the closure of $\mathbb{H}^\alpha$ in $L^2$, it is easy to check that the support is in fact $L^2$ itself. Furthermore, the small ball probabilities $\varphi^B$ for this process have a well-known behavior, that is $-\log \mathbf{P}(\|X_\alpha\|_2 < \varepsilon)$ is of the order of $\varepsilon^{-1/\alpha}$ as $\varepsilon \to 0$, see for instance [8], Theorem 4.

*Gaussian white noise model.* To simplify the formulation of the upper-bound results, we will assume that we are in a particularly simple model, namely the Gaussian white model. In this model the data $X^{(n)}$ is given by

$$dX^{(n)}(t) = f(t)dt + \frac{1}{\sqrt{n}}dW(t), \quad t \in [0,1], \tag{6}$$

for some $f$ in $L^2[0,1]$ and $W$ standard Brownian motion. Let us denote, for any positive real numbers $\alpha$ and $\beta$,

$$r_n^{\alpha,\beta} \triangleq n^{-\frac{\alpha \wedge \beta}{2\alpha+1}}. \tag{7}$$

In the sequel the notation $\lesssim$ is used for "smaller than or equal to a universal constant" and $\gtrsim$ is defined similarly.

**Theorem 2.** *Let $\beta > 0$, $L > 0$ and suppose the data is generated according to (6). Let the prior process be defined by (5) with $\alpha > 0$. Let $f_0$ be in $\mathcal{F}_{\beta,L}$ and let the rate $r_n^{\alpha,\beta}$ be defined by (7). Let $\varepsilon_n$ and $\zeta_n$ be such that*

$$\varphi_{f_0}(\varepsilon_n) \leq n\varepsilon_n^2 \quad \text{and} \quad \zeta_n \leq \varphi_{f_0}^{-1}(9n\varepsilon_n^2).$$

*Then for $M$ large enough,*

$$\mathbf{E}_0 \Pi(\zeta_n \leq \|f - f_0\|_2 \leq M\varepsilon_n \mid X^{(n)}) \to 1,$$

*as $n \to +\infty$. For any $f_0$ in $\mathcal{F}_{\beta,L}$, one can choose $\varepsilon_n$ such that $\varepsilon_n \lesssim r_n^{\alpha,\beta}$ and, if $\alpha \leq \beta$, one can choose $\zeta_n$ such that $\zeta_n \gtrsim r_n^{\alpha,\beta}$. Furthermore, if $\beta < \alpha$, there exists $f_0$ in $\mathcal{F}_{\beta,L}$ such that, for $p > 1 + \beta/2$ and $M$ large enough, as $n \to +\infty$,*

$$\mathbf{E}_0 \Pi(r_n^{\alpha,\beta} \log^{-p} n \leq \|f - f_0\|_2 \leq Mr_n^{\alpha,\beta} \mid X^{(n)}) \to 1.$$

The first convergence result is essentially a consequence of Theorem 3.4 in [14] for the upper-bound and of Theorem 1 for the lower bound. The second part of the statement reveals that there are indeed functions in the class such that the posterior rate is $r_n^{\alpha,\beta}$, up to a log-factor if $\beta < \alpha$. In this sense the rate can be said to be *optimal* (up to a log-factor) over $\mathcal{F}_{\beta,L}$.

It is interesting to compare these results to the ones obtained by [17] and [2], where the authors also study estimation in model (6) from the Bayesian perspective. Both works obtain the upper-bound result on $\varepsilon_n$ for priors defined by (5) by different methods but they do not consider the question of optimality of the rate $r_n^{\alpha,\beta}$ when $\alpha \neq \beta$. In [2], the focus is on the question of adaptation when one puts also a prior on $\alpha$ and the authors obtain the minimax rate for the resulting prior for unknown $\beta$ under some conditions. In [17], an interesting observation about non-optimality is made, but in a rather different direction than ours, the author noting that though the prior (5) leads to the minimax rate for $\alpha = \beta$, both the prior and the posterior put mass zero on the Sobolev space $\{f = (f_k)_{k \geq 1}, \sum_{k \geq 1} k^{2\beta} f_k^2 < +\infty\}$ the true function $f_0$ belongs to.

*Remark 1.* If $\alpha \leq \beta$, the precise rate of convergence of the posterior is, up to constants, equal to $r_n^{\alpha,\alpha} = n^{-\alpha/(2\alpha+1)}$. If $\alpha > \beta$, more information on $f_0$ (for instance about the rate of decrease of its Fourier coefficients) is needed to evaluate the RKHS-approximation term and eventually obtain an explicit expression of the rate, see for example the special "worst-case" function $f_0$ considered in the proof of the theorem.

*Remark 2.* It is natural to ask whether it is possible to avoid the log-factor for the lower bound. The answer is yes if one allows sequences of functions: it can be checked that there exists a sequence $f_{0,n}$ in $\mathcal{F}_{\beta,L}$, where the function $f_{0,n}$ has only one properly chosen non-zero Fourier coefficient, such that, for $M$ large enough, $\mathbf{E}_{f_{0,n}} \Pi(r_n^{\alpha,\beta}/M \leq \|f - f_{0,n}\|_2 \mid X^{(n)})$ tends to 1 as $n \to +\infty$.

*Proof of Theorem 2.* The fact that the posterior concentrates in a ball of radius $M\varepsilon_n$ for the $\|\cdot\|_2$-norm is the conclusion of Theorem 3.4 in [14]. The explicit upper-bound for $\varepsilon_n$ is obtained as follows. Denoting $f_K = \sum_{k=1}^{K} f_{0,k} e_k(\cdot)$, note that $f_K$ belongs to $\mathbb{H}^\alpha$. Since $f_0$ belongs to $\mathcal{F}_{\beta,L}$, it holds

$$\|f_K - f_0\|_2^2 = \sum_{p \geq K+1} f_{0,p}^2 \leq K^{-2\beta} \sum_{p \geq K+1} p^{2\beta} f_{0,p}^2 \leq L^2 K^{-2\beta}$$

$$\|f_K\|_{\mathbb{H}^\alpha}^2 = \sum_{p=1}^{K} p^{1+2\alpha} f_{p,0}^2 \leq K^{(1+2\alpha-2\beta)\vee 0} \sum_{p=1}^{K} p^{2\beta} f_{p,0}^2 \leq L^2 K^{(1+2\alpha-2\beta)\vee 0}.$$

Let us now choose $K = \varepsilon_n^{-1/\beta}$. The last display then implies that the approximation part $\varphi_{f_0}^A(\varepsilon_n)$ of the concentration function is at most $\varepsilon_n^{-(1+2\alpha-2\beta)/\beta \wedge 0}$. On the other hand, the small ball probability $\varphi^B(\varepsilon_n)$ is at most constant times $\varepsilon_n^{-1/\alpha}$ for $n$ large enough as noted at the beginning of this Section. Hence

$$\varphi_{f_0}(\varepsilon_n) \lesssim \varepsilon_n^{-1/\alpha} + \varepsilon_n^{-(1+2\alpha-2\beta)/\beta \wedge 0}.$$

If we choose $n\varepsilon_n^2$ equal to the latter quantity we get $\varepsilon_n \lesssim n^{-\alpha \wedge \beta/(2\alpha+1)} = r_n^{\alpha,\beta}$.

To obtain the lower bound result, we apply Theorem 1. Simple calculations reveal that for model (6), the set $B_{KL}(f_0, \varepsilon)$ coincides with $\{f, \|f - f_0\|_2 < \varepsilon\}$, see Lemma 6 in [6] and thus $\Pi(B_{KL}(f_0, 2\varepsilon_n)) = \Pi(\|f - f_0\|_2 \leq 2\varepsilon_n)$. Now apply the remark after Theorem 1 to obtain that if $\varphi_{f_0}(\varepsilon_n) \leq n\varepsilon_n^2$, then any $\zeta_n$ such that $\varphi_{f_0}(\zeta_n) \geq 9n\varepsilon_n^2$ is a lower bound for the rate. To obtain a more explicit form for $\zeta_n$, we distinguish the cases $\alpha \leq \beta$ and $\beta \leq \alpha$.

In the case $\alpha \leq \beta$, let us use the fact that

$$\varphi_{f_0}(\zeta_n) \geq -\log \Pi(\|f\|_2 < \zeta_n) \gtrsim \zeta_n^{-1/\alpha},$$

where the last inequality is obtained using the asymptotics of the small ball probability of $X_\alpha$. Thus the condition $\varphi_{f_0}(\zeta_n) \geq 9n\varepsilon_n^2$ is satisfied if $\zeta_n$ is equal to constant times $n^{-\alpha/(2\alpha+1)} = r_n^{\alpha,\alpha}$, since $\varepsilon_n$ can be chosen equal to constant times $r_n^{\alpha,\beta} = r_n^{\alpha,\alpha}$.

In the case $\alpha > \beta$, let us define $f_0$ by specifying its Fourier coefficients as

$$f_{0,k}^{-1} = k^{1/2+\beta}(1 + \log k)^{1/2} \log \log k, \qquad (k \geq 1).$$

Note that the series $\sum k^{2\beta} f_{0,k}^2$ converges so without loss of generality one can assume that $f_0$ belongs to $\mathcal{F}_{\beta,L}$ (otherwise consider $af_0$ for $a > 0$ small enough). Moreover, one just needs to prove the lower bound result, the upper-bound resulting from what precedes. In the remainder of the proof the rate $\varepsilon_n$ is thus taken equal to $Cr_n^{\alpha,\beta}$ for some constant $C > 0$.

Let us denote $\zeta_n = \delta_n \varepsilon_n$, where $\delta_n \to 0$ is to be chosen, and let us bound from below $\varphi_{f_0}(\zeta_n)$. We have $\varphi_{f_0}(\zeta_n) \geq \varphi_{f_0}^A(\zeta_n)$. Let $h$ be in the RKHS $\mathbb{H}^\alpha$ of the prior with $\|h - f_0\|_2 < \zeta_n$. Then, for any $k(n) \geq 1$, using the inequality $(x+y)^2 \geq x^2/2 - y^2$ valid for all reals $x$ and $y$,

$$\begin{aligned} \|h\|_{\mathbb{H}}^2 &= \sum_{k \geq 1} k^{1+2\alpha} h_k^2 \geq \sum_{k=1}^{k(n)} k^{1+2\alpha}(h_k - f_{0,k} + f_{0,k})^2 \\ &\geq \frac{1}{2} \sum_{k=1}^{k(n)} k^{1+2\alpha} f_{0,k}^2 - \sum_{k=1}^{k(n)} k^{1+2\alpha}(h_k - f_{0,k})^2. \end{aligned}$$

That is, with the notation $S(K) = \sum_{k=1}^K k^{1+2\alpha} f_{0,k}^2$, using that $\|h - f_0\|_2 < \zeta_n$,

$$\|h\|_{\mathbb{H}}^2 \geq \frac{1}{2} S(k(n)) - k(n)^{1+2\alpha} \zeta_n^2. \qquad (8)$$

Let us choose $k(n) = n^{1/(1+2\alpha)} \log n$ and $\delta_n = \log^{-p} n$ for some $p > 0$. Using the explicit form of the $f_{0,k}$'s, one obtains, denoting $l_n = \log \log n$, that $S(k(n)) \gtrsim k(n)^{1+2\alpha-2\beta} l_n^{-2} \log^{-1} n$. Thus

$$\begin{aligned} S(k(n)) &\gtrsim n\varepsilon_n^2 l_n^{-2} \log^{2\alpha-2\beta} n \\ k(n)^{1+2\alpha} \zeta_n^2 &= n\varepsilon_n^2 \log^{2\alpha+1-2p} n. \end{aligned}$$

Since $\alpha > \beta$, the first of these two terms is of larger order than $n\varepsilon_n^2$. As soon as $2p > 1 + 2\beta$, it is also of larger order than the last term in the preceding display. Minimizing Equation (8) in $h$, we conclude that in this case,

$$\varphi_{f_0}^A(\zeta_n) \gtrsim n\varepsilon_n^2 l_n^{-2} \log^{2\alpha-2\beta} n.$$

Thus $\varphi_{f_0}(\zeta_n)$ divided by $n\varepsilon_n^2$ tends to infinity. In view of the Remark after Theorem 1, we obtain that $\zeta_n = \delta_n \varepsilon_n$ is a lower bound for the rate, which concludes the proof. □

### 3.2. *The $\mathcal{C}^0[0,1]$-setting and Riemann-Liouville type priors*

In this subsection we obtain new upper and lower bounds for posterior rates in the following model of density estimation. The observations $X_1, \ldots, X_n$ are a random sample from a *positive* density $f_0$ on the interval $[0,1]$. Let us denote $w_0 = \log f_0$, so that $f_0 = e^{w_0}$.

Now let us explain how we construct a prior on positive densities $f$, following the approach considered in [14]. To any continuous function $w$ on the interval $[0,1]$, we associate the density $p_w$ (that is a nonnegative function which integrates to 1) defined by

$$p_w(s) = \frac{e^{w(s)}}{\int_0^1 e^{w(u)}du}, \quad s \in [0,1].$$

Let $W$ be a Gaussian process defining a prior $\Pi_w$ on $\mathcal{C}^0[0,1]$. Then the quantity $p_W$ defines a random (non-Gaussian) density. The corresponding prior on the set of densities is denoted by $\Pi_{p_w}$. As Gaussian prior $W$ we choose the process $X_t^\alpha$ defined below.

First, let us define the Riemann-Liouville process of parameter $\alpha > 0$ as

$$R_t^\alpha = \int_0^t (t-s)^{\alpha - 1/2} dB(s), \quad t \in [0,1], \tag{9}$$

where $B$ is standard Brownian motion. Then the process prior, which we call the *Riemann-Liouville type process* (RL-type process), is defined as

$$X_t^\alpha = R_t^\alpha + \sum_{k=0}^{\underline{\alpha}+1} Z_k t^k, \quad t \in [0,1],$$

where $Z_0, \ldots, Z_{\underline{\alpha}+1}, R_t$ are independent, $Z_i$ is standard normal and $R_t^\alpha$ is the Riemann-Liouville process of parameter $\alpha$. Note that if $\alpha = 1/2$ then $R_t^\alpha$ is simply standard Brownian motion and if $\{\alpha\} = 1/2$, with $\{\alpha\} \in [0,1)$ the integer part of $\alpha$, then $R_t^\alpha$ is a $k$-fold integrated Brownian motion. It can be checked that the support in $\mathcal{C}^0[0,1]$ of $X_t^\alpha$ is the whole space $\mathcal{C}^0[0,1]$ (it is in fact in order to get the whole space as support that one adds the polynomial part), see [14], Section 4 and Theorem 4.3.

Let us denote by $\varphi_{w_0}$ the concentration function associated to the process $X_t^\alpha$ and the continuous function $w_0$. Upper-bounds on $\varphi_{w_0}$ used in the proof of the next Theorem to get explicit upper bound rates are obtained in Section 4.1.

**Theorem 3.** *Suppose that $w_0 = \log f_0$ belongs to the Hölder class $\mathcal{C}^\beta[0,1]$ for some $\beta > 0$ and let the prior on densities be the distribution $\Pi_{p_w}$ of $p_{X^\alpha}$, where $X^\alpha$ is a Riemann-Liouville type process of parameter $\alpha > 0$. Then there exist finite constants $C_1, C_2 > 0$ such that, if $\varepsilon_n$ and $\zeta_n$ are such that*

$$\varphi_{w_0}(\varepsilon_n) \leq n\varepsilon_n^2 \quad \text{and} \quad \zeta_n \leq C_1 \varphi_{w_0}^{-1}(C_2 n\varepsilon_n^2),$$

*then for M large enough, as $n \to +\infty$,*

$$\mathbf{E}_0 \Pi_{p_w}(h(f, f_0) \le M\varepsilon_n \mid X^{(n)}) \to 1,$$
$$\mathbf{E}_0 \Pi_{p_w}(\|f - f_0\|_\infty \ge \zeta_n \mid X^{(n)}) \to 1,$$

*where h is Hellinger's distance. Moreover, one can choose $\varepsilon_n$ such that $\varepsilon_n \lesssim r_n^{\alpha,\beta}$ if $\{\alpha\} = 1/2$ or $\alpha$ does not belong to $\beta + 1/2 + \mathbb{N}$ and $\varepsilon_n \lesssim n^{-\beta/(2\alpha+1)} \log n$ otherwise.*

These results describe in a rather complete way the rate of convergence of the posterior for the prior $\Pi_{p_w}$ constructed from the Riemann-Liouville prior, for all values of the parameter $\alpha$ in $(0, +\infty)$. Also, from the upper-bounds point of view, it improves on Theorem 4.3 in [14], where $\alpha = \beta$ is needed.

Note that, while upper-bound rates are obtained for Hellinger's distance, the lower bounds are in terms of the uniform norm. To obtain the lower bounds, the uniform norm is in a way the most natural (and easiest) distance to work with since it is the norm on the Banach space where the prior lives and in the proof, the idea will be indeed to apply Lemma 1 with sets of the form $B_n = \{f, \|f - f_0\|_\infty \le \zeta_n\}$. For upper-bounds, Hellinger's distance is a rather natural choice since it is a natural testing distance for i.i.d. observations in view of the theory of [5]. A natural extension of our results would be to obtain results in terms of a common distance on the parameter space. While such a refinement is beyond the scope of the present contribution, we hope that future papers will answer this type of question.

In the above Theorem, explicit bounds for $\varepsilon_n$ are obtained using explicit upper-bounds for the concentration functions obtained in Section 4.1. It should also be possible to obtain explicit bounds for $\zeta_n$ in the spirit of those of Theorem 2 by bounding the concentration function from below. One difficulty here with respect to Theorem 2 is the presence of the extra polynomial part in the definition of the process, which makes the evaluations even for the small ball term more difficult. We will not further discuss this issue here but note only that in some simple cases, an explicit expression for $\zeta_n$ follows quite directly from what precedes.

*Remark 3.* For Brownian motion released at zero $X_t = B_t + Z_0$, by a slight adaptation of the preceding, one can obtain an explicit evaluation of $\zeta_n$ and show that if $f_0$ is smooth enough, more precisely if $\beta \ge 1/2$, there exist a constant $m$ such that, as $n \to +\infty$,

$$\Pi_{p_w}(\|f - f_0\|_\infty \ge mn^{-1/4}|X^{(n)}) \to 1.$$

Note that $X_t$ is almost the RL-type process with $\alpha = 1/2$ except for the term $tZ_1$. It can still be checked that the support in $\mathcal{C}^0[0,1]$ of this process is the full space $\mathcal{C}^0[0,1]$, see [14], Theorem 4.1, and that Theorem 3 still holds, following the same proof. We always have $\varphi_{f_0}(\varepsilon) \ge \varphi^B(\varepsilon) = -\log \mathbf{P}(\|X\|_\infty < \varepsilon)$. But on the event that $\|X\|_\infty < \varepsilon$, we have $|X_0| = |Z_0| < \varepsilon$ thus $\|B\|_\infty < 2\varepsilon$.

Using the behavior of the small ball probability of Brownian motion, we obtain that there exists a constant $C$ such that $\varphi_{f_0}(\varepsilon) \geq C\varepsilon^{-2}$ for $\varepsilon$ small enough and thus $\varphi_{f_0}^{-1}(u) \gtrsim u^{-1/2}$. Hence $\zeta_n$ can be chosen equal to constant times $(n\varepsilon_n^2)^{-1/2}$, where $\varepsilon_n$ is the upper-bound rate obtained for $X_t$. By the smoothness assumption on $f_0$, the rate $\varepsilon_n$ can be chosen equal to a constant times $r_n^{1/2,1/2} = n^{-1/4}$, which yields the announced result on $\zeta_n$.

The following Lemmas are used in the proof of Theorem 3. The proof of Lemma 4 can be found in Section 4.2.

**Lemma 4.** *Let $\varphi_w$ denote the concentration function associated to the process $X_t^\alpha$ and the function $w \in \mathcal{C}^0[0,1]$ and let $\rho$ denote both the real $\rho$ and the constant function equal to $\rho$. Then for any $\varepsilon > 0$,*

$$\varphi_{w_0+\rho}(\varepsilon) \geq \varphi_{w_0}(\varepsilon) + \rho^2 - 2(|w_0(0)| + \varepsilon)|\rho|.$$

**Lemma 5** (Lemma 3.1 in [14])**.** *For any $v, w$ elements of $\mathcal{C}^0[0,1]$,*

$$h(p_v, p_w) \leq \|v - w\|_\infty e^{\|v-w\|_\infty/2}$$
$$K(p_v, p_w) \vee V(p_v, p_w) \lesssim \|v - w\|_\infty^2 e^{\|v-w\|_\infty/2}(1 + \|v - w\|_\infty)^2.$$

*Proof of Theorem 3.* The fact that any $\varepsilon_n$ such that $\varphi_{w_0}(\varepsilon_n) \leq n\varepsilon_n^2$ is an upper bound for the rate is Theorem 3.1 in [14]. Now Theorem 4 in the Appendix enables to get the explicit expression of $\varepsilon_n$ in terms of $r_n^{\alpha,\beta}$.

To obtain the lower bound result, we show that if $\zeta_n = C_3^{-1} \varphi_{w_0}^{-1}(C_1 n\varepsilon_n^2)$ for large enough constants $C_1$ and $C_3$ then $\Pi_{p_w}(\|f - f_0\|_\infty \leq \zeta_n) \leq \exp(-3n\varepsilon_n^2)$. This is enough to obtain the lower bound statement, since then one can apply Lemma 1 with $B_n = \{f, \|f - f_0\|_\infty \leq \zeta_n\}$. The prior probability on the Kullback-Leibler type neighborhood is bounded from below using Lemma 5 to obtain a neighborhood in terms of the $\|\cdot\|_\infty$-norm, and finally, due the fact that the support of $X^\alpha$ in $\mathcal{C}^0[0,1]$ is the whole space $\mathcal{C}^0[0,1]$, Lemma 2 can be used.

Let $A_n$ be the set $\{w \in \mathcal{C}^0[0,1], \|p_w - f_0\|_\infty \leq \zeta_n\}$. Since $\zeta_n \to 0$ and $f_0 \geq \rho > 0$ for some $\rho > 0$, it holds $2\|f_0\|_\infty \geq p_w \geq \rho/2 > 0$ on $A_n$ for $n$ large enough. Since the logarithm is a Lipshitz function on the interval $[\rho/2, 2\|f_0\|_\infty]$, one gets, on $A_n$, for some $d > 0$,

$$\|\log p_w - \log f_0\|_\infty \leq d\|p_w - f_0\|_\infty \leq d\zeta_n.$$

Noting that

$$\|\log p_w - \log f_0\|_\infty = \|\log \frac{e^w}{\int e^w} - w_0\|_\infty = \|w - w_0 - \log \int e^w\|_\infty,$$

one obtains that, on $A_n$, it holds $\|w - w_0 - Z_w\|_\infty \leq d\zeta_n$, where $Z_w$ is a constant function and $|Z_w| \leq \|w\|_\infty$. We shall use the fact that with high probability, this value is not too large. Note that if $w$ is in $A_n$ and $\|w\|_\infty \leq C\sqrt{n}\varepsilon_n$ then $w$

belongs to $\cup_{k=-N}^{N} B_k$, where $B_k = \{w, \|w - w_0 - c_k\|_\infty \leq 2d\zeta_n\}$ with $c_k = kd\zeta_n$ and $N$ the smallest integer larger than $C\sqrt{n}\varepsilon_n/(d\zeta_n)$. Thus

$$\Pi_{p_w}(\|f - f_0\|_\infty \leq \zeta_n) = \Pi_w(\|p_w - f_0\|_\infty \leq \zeta_n)$$
$$\leq \sum_{k=-N}^{N} \Pi_w(\|w - w_0 - c_k\|_\infty \leq 2d\zeta_n) + \Pi_w(\|w\|_\infty > C\sqrt{n}\varepsilon_n).$$

It is an easy consequence of Borell's inequality, see [3] or [16], Proposition A.2.1, that $\Pi_w(\|w\|_\infty > C\sqrt{n}\varepsilon_n)$ is bounded above by $\exp(-4n\varepsilon_n^2)$ for $C$ large enough. Now due to Lemma 2,

$$\Pi_w(\|w - w_0 - c_k\|_\infty \leq 2d\zeta_n) \leq \exp(-\varphi_{w_0+c_k}(2d\zeta_n)).$$

Let $I_1$ be the set of indexes $k$ such that $|c_k| \leq 4|w_0(0)|$ and $I_2$ the set of indexes such that $|c_k| > 4|w_0(0)|$. According to Lemma 4, we have for $n$ large enough

$$\varphi_{w_0+c_k}(2d\zeta_n) \geq \begin{cases} \varphi_{w_0}(2d\zeta_n) - 9|w_0(0)|^2, & \text{if } k \in I_1, \\ \varphi_{w_0}(2d\zeta_n) + c_k^2/2, & \text{if } k \in I_2. \end{cases}$$

Thus for some $C_4 > 0$, it holds

$$\Pi_{p_w}(\|f - f_0\|_\infty \leq \zeta_n)$$
$$\lesssim \zeta_n^{-1} e^{9|w_0(0)|^2 - \varphi_{w_0}(2d\zeta_n)} + \sum_{\zeta_n^{-1} \lesssim |k| \lesssim N} e^{-k^2 d^2 \zeta_n^2/2 - \varphi_{w_0}(2d\zeta_n)} + e^{-4n\varepsilon_n^2}$$
$$\leq C_4(\zeta_n^{-1} e^{-\varphi_{w_0}(2d\zeta_n)} + e^{-4n\varepsilon_n^2}).$$

Using the behavior of the small ball probability for the process at stake, we have that $\varphi_{w_0}(2d\zeta_n) \gtrsim \zeta_n^{-1/\alpha}$ hence for $n$ large enough it holds $\varphi_{w_0}(2d\zeta_n) + 2\log \zeta_n \geq \varphi_{w_0}(2d\zeta_n)/2$. Thus the last display is bounded from above by $2C_4 \exp(-4n\varepsilon_n^2)$ as soon as $\varphi_{w_0}(2d\zeta_n) \geq 8n\varepsilon_n^2$, which concludes the proof. □

## 4. Appendix

### *4.1. Concentration function of RL-type processes: upper bounds*

In this subsection, we establish an upper-bound result on the concentration function of the RL-type process which is of independent interest and which is used in the proof of Theorem 3 to get explicit upper bound rates.

First let us introduce the classical notion of fractional integral, whose definition is as follows. For $\alpha > 0$ and $f$ a continuous function on $[0, 1]$, the fractional integral of order $\alpha$ is defined as

$$I_{0+}^\alpha f(t) = \int_0^t (t-s)^{\alpha-1} f(s) ds,$$

for any $t$ in $[0,1]$. If $t < 0$, we set $I^\alpha_{0+} f(t) = 0$.

We shall use the following two Lemmas, the second enabling to handle a case discarded by the first one (namely the case where $\alpha + \lambda = 1$). In the next statement, the symbol $*$ stands for the usual convolution between functions.

**Lemma 6** (Lemma 5.2 in [14]). *Let $\lambda \in [0,1]$ and $\alpha \in [0,1)$ be such that $\alpha + \lambda \in (0,2)$ and $\alpha + \lambda \neq 1$. If $f \in \mathcal{C}^\lambda[0,1]$ and $g \in L_1(\mathbb{R})$ has compact support and satisfies $\int g(u)du = 0$ and, in the case that $\alpha + \lambda > 1$, also $\int ug(u)du$, then*

$$\|I^\alpha_{0+}(f*g)\|_\infty \lesssim \int |u|^{\alpha+\lambda} |g(u)| du.$$

**Lemma 7.** *Let $\delta \in (0,1)$ and $f \in \mathcal{C}^\delta[0,1]$. If $g \in L_1(\mathbb{R})$ has compact support and satisfies $\int g(u)du = 0$ then*

$$\|I^{1-\delta}_{0+}(f*g)\|_2^2 \lesssim \int u^2 \{1 + \log^2(1+|u|^{-1})\} g(u)^2 du.$$

The proof of Lemma 7 can be found in Section 4.2.

**Theorem 4.** *Suppose $f_0$ belongs to $\mathcal{C}^\beta[0,1]$, with $\beta > 0$. The concentration function $\varphi_{f_0}$ associated to the process $X^\alpha_t$ satisfies, if $0 < \alpha \leq \beta$, that $\varphi_{f_0}(\varepsilon) = O(\varepsilon^{-1/\alpha})$ as $\varepsilon \to 0$. In the case that $\alpha > \beta$, as $\varepsilon \to 0$,*

$$\varphi_{f_0}(\varepsilon) = \begin{cases} O(\varepsilon^{-\frac{2\alpha-2\beta+1}{\beta}}) & \text{if } \{\alpha\} = 1/2 \text{ or } \alpha \notin \beta + \frac{1}{2} + \mathbb{N}, \\ O(\varepsilon^{-\frac{2\alpha-2\beta+1}{\beta}} \log(1/\varepsilon)) & \text{otherwise.} \end{cases}$$

This extends Theorem 4.3 in [14] in the case that $\alpha \neq \beta$. There is an extra difficulty in the case where $\alpha - \beta - 1/2$ is an integer and $\{\alpha\}$ is not $1/2$, resulting in the presence of the extra log-factor. Roughly, the difficulty arises from the fact that, if $\alpha \in (0,1)$ and $\lambda \in [0,1]$, the fractional integral $I^\alpha_{0+}$ does map $\mathcal{C}^\lambda[0,1] \to \mathcal{C}^{\lambda+\alpha}[0,1]$ only if $\alpha + \lambda \neq 1$, see [7]. Lemma 7 enables us to deal with the case where $\alpha + \lambda = 1$ is an integer.

*Proof of Theorem 4.* Let us denote by $Z = X^\alpha - R^\alpha$ the polynomial part of $X^\alpha$ and by $\mathbb{H}^\alpha$ the RKHS of $R^\alpha$. The proof is quite similar to the one of Theorem 4.3 in [14] and the starting point is identical. Using Theorem 2.3 in [14], the initial step of the proof is to bound from above the concentration function $\varphi_{f_0}(2\varepsilon)$ by a multiple of the sum $\varphi_{f_0-P}(\varepsilon/2, R^\alpha) + \varphi_P(\varepsilon/2, Z)$ with the polynomial $P$ to be chosen in the RKHS $\mathbb{H}^Z$ of $Z$. The spaces $\mathbb{H}^Z$ and $\mathbb{H}^\alpha$ are known explicitly. The space $\mathbb{H}^Z$ is the set of polynomials $P_\xi = \sum_{i=0}^{\alpha+1} \xi_i t^i$ equipped with the norm $\|P_\xi\|^2_{\mathbb{H}^Z} = \sum_{i=0}^{\alpha+1} \xi_i^2$. The RKHS $\mathbb{H}^\alpha$ is the space $I^{\alpha+1/2}_{0+}(L^2[0,1])$ with associated norm $\|I^{\alpha+1/2}_{0+} f\|_{\mathbb{H}^\alpha} = \|f\|_2/\Gamma(\alpha+1/2)$, where $\Gamma$ is the Gamma function, due to Theorem 4.2 in [14].

Let us check that for the process $X^\alpha$, the small ball term $\varphi^B(\varepsilon)$ is bounded above by a constant times $\varepsilon^{-1/\alpha}$ for $\varepsilon$ small enough. Indeed, it is known that for the Riemann-Liouville process $R^\alpha$, the quantity $-\log \mathbf{P}(\|R^\alpha\|_\infty < \varepsilon)$ behaves as a constant times $\varepsilon^{-1/\alpha}$ as $\varepsilon \to 0$, see [10]. Moreover, for any integer $k$, the quantity $\mathbf{P}(\|Z_k t^k\|_\infty < \varepsilon)$ behaves as a constant times $\varepsilon$, which by independence of $Z$ and $R^\alpha$ implies that $-\log \mathbf{P}(\|X^\alpha\|_\infty < \varepsilon)$ is smaller than a constant times $\varepsilon^{-1/\alpha}$ for $\varepsilon$ small enough.

Now we study the RKHS-approximation term $\varphi_{f_0}^A(\varepsilon)$. There are different cases depending on the value of $\{\alpha\}$ which are: $\{\alpha\} \in (0, 1/2)$, $\{\alpha\} \in (1/2, 1)$, $\{\alpha\} = 0$ and $\{\alpha\} = 1/2$. We will focus on the first case, that is $\{\alpha\} \in (0, 1/2)$, the other cases being similar. Also, we assume that $\alpha \geq \beta$, the case $\alpha < \beta$ being similar though easier since the small ball term dominates in that case.

Thus we focus on the RKHS-approximation term $\varphi_{f_0}^A(\varepsilon)$ in the case where $\alpha \geq \beta$ and $\{\alpha\} \in (0, 1/2)$. Let $\phi$ be a smooth, compactly supported kernel, of sufficiently large order and for $\sigma > 0$ define $\phi_\sigma(t) = \sigma^{-1}\phi(t/\sigma)$. Since $f_0 \in \mathcal{C}^\beta$, we have $\|f_0 - f_0 * \phi_\sigma\|_\infty \lesssim \sigma^\beta$ thus $\|\{f_0 - P\} - \{f_0 * \phi_\sigma - P\}\|_\infty \leq \varepsilon$ if $\sigma = C\varepsilon^{1/\beta}$ for some constant $C$.

Let us write Taylor's theorem in the form

$$f_0 * \phi_\sigma(t) = \sum_{k=0}^{\underline{\alpha}} \frac{(f_0 * \phi_\sigma)^{(k)}(0)}{k!} t^k + I_{0+}^{\alpha+1/2} I_{0+}^{1/2-\{\alpha\}}(f_0^{(\underline{\beta})} * \phi_\sigma^{(\underline{\alpha}-\underline{\beta}+1)}).$$

For the polynomial $P$ let us choose the polynomial part in the preceding display. Its squared RKHS-norm $\|P\|_{\mathbb{H}^Z}^2$ is proportional to $\sum_{k=0}^{\underline{\alpha}} (f_0 * \phi_\sigma)^{(k)}(0)^2$. The term of largest order is $(f_0 * \phi_\sigma)^{(\underline{\alpha})}(0)^2 = f_0^{(\underline{\beta})} * \phi_\sigma^{(\underline{\alpha}-\underline{\beta})}(0)^2$. Note that, since $f_0$ is in $\mathcal{C}^\beta$, denoting by $\{\beta\}$ the fractional part of $\beta$,

$$|f_0^{(\underline{\beta})} * \phi_\sigma^{(\underline{\alpha}-\underline{\beta})}(0)| = |\int \{f_0^{(\underline{\beta})}(0-s) - f_0^{(\underline{\beta})}(0)\} \phi_\sigma^{(\underline{\alpha}-\underline{\beta})}(s) ds|$$
$$\leq \int |s|^{\{\beta\}} \phi_\sigma^{(\underline{\alpha}-\underline{\beta})}(s) ds \lesssim \sigma^{\beta-\underline{\alpha}}.$$

Hence $\|P\|_{\mathbb{H}^Z}^2 \lesssim \sigma^{2\beta-2\underline{\alpha}} \leq \sigma^{-1-2\alpha+2\beta}$. Now notice that $f_0 * \phi_\sigma - P$ belongs to $\mathbb{H}^\alpha$ and has RKHS-norm proportional to $\|I_{0+}^{1/2-\{\alpha\}}(f_0^{(\underline{\beta})} * \phi_\sigma^{(\underline{\alpha}-\underline{\beta}+1)})\|_2$. Thus, in the case where $1/2 - \{\alpha\} + \{\beta\} \neq 1$, we can use Lemma 6 to get

$$\|I_{0+}^{1/2-\{\alpha\}}(f_0^{(\underline{\beta})} * \phi_\sigma^{(\underline{\alpha}-\underline{\beta}+1)})\|_2$$
$$\lesssim \int |u|^{1/2-\{\alpha\}+\{\beta\}} \sigma^{-\underline{\alpha}+\underline{\beta}-2} |\phi^{(\underline{\alpha}-\underline{\beta}+1)}(u/\sigma)| du \lesssim \sigma^{-1/2-\alpha+\beta}.$$

Thus $\|f_0 * \phi_\sigma - P\|_{\mathbb{H}^\alpha}^2 \lesssim \sigma^{-1-2\alpha+2\beta} \lesssim \varepsilon^{-(2\alpha-2\beta+1)/\beta}$. The small ball term $\varphi^B$ being of smaller order, the concentration function is at most of order $\varepsilon^{-(2\alpha-2\beta+1)/\beta}$, which concludes the proof in this case.

If $1/2 - \{\alpha\} + \{\beta\} = 1$, let us apply Lemma 7 to obtain

$$\|I_{0+}^{1/2-\{\alpha\}}(f_0^{(\underline{\beta})} * \phi_\sigma^{(\underline{\alpha}-\underline{\beta}+1)})\|_2^2 \lesssim \sigma^{-2\underline{\alpha}+2\underline{\beta}-1} \int v^2 \phi^{(\underline{\alpha}-\underline{\beta}+1)}(v)\{1+\log^2(1+\frac{1}{|\sigma v|})\} dv.$$

Using the inequality $1 + (|\sigma v|)^{-1} \le \sigma^{-1}(1 + |v|^{-1})$ valid for $0 < \sigma < 1$, one obtains that the norm of $f_0 * \phi_\sigma - P$ in $\mathbb{H}^\alpha$ is bounded by constant times $\log(\sigma^{-1})\sigma^{-1/2-\alpha+\beta}$, which concludes the proof. □

### 4.2. Proof of the Lemmas

*Proof of Lemma 3.* The concentration function $\varphi_{f_0}$ is the sum $\varphi_{f_0}^A + \varphi^B$ of two decreasing functions, see Equation (3). Let us show that $\varphi^B$ is strictly decreasing that is $\varphi^B(\varepsilon) > \varphi^B(\varepsilon')$ if $\varepsilon' > \varepsilon$. It suffices to see that $\mathcal{C} = \{\gamma \in \mathbb{B},\ \varepsilon < \|\gamma\| < \varepsilon'\}$ receives positive mass under the law of $Z$. Since $Z$ is non-degenerate, its RKHS $\mathbb{H}$ contains a non-zero element $h_1$. For some $\lambda > 0$, the element $\lambda h_1 \in \mathbb{H}$ belongs to the set $\mathcal{C}$. Since $\mathcal{C}$ is an open set, there exists $\eta > 0$ such that the ball $B(\lambda h_1, \eta)$ is included in $\mathcal{C}$. Thus to obtain the strict monotonicity, it suffices to show that the probability $\mathbf{P}(\|Z - h\| < \eta)$ of an arbitrary open ball centered around an element of $h \in \mathbb{H}$ is positive.

This can be proved as follows. Using Cameron-Martin change of variables formula, see e.g. Lemma 3.2 in [15], it suffices to prove that any ball centered at zero of positive radius receives positive mass under the law of $Z$. Since $\mathbb{B}$ is separable, for any $r > 0$ the space $\mathbb{B}$ is the union of a countable number of balls of radius $r$ thus at least one of these balls say $B(x, r)$ receives positive mass. Let $Z'$ be an independent copy of the process $Z$, we have

$$0 < \mathbf{P}(Z \in B(x,r))\mathbf{P}(Z' \in B(x,r)) \le \mathbf{P}(Z - Z' \in B(0, 2r)).$$

Now note that $(Z - Z')/\sqrt{2}$ has the same distribution as $Z$ (these are Gaussian processes with the same covariance function) thus $\mathbf{P}(Z \in B(0, r))$ is positive for any $r > 0$.

Now just the convexity statement remains to prove. Using the fact that the function $h \to \|h\|_\mathbb{H}^2$ is convex together with the definition of the infimum, one gets that $\varphi_{f_0}^A$ is convex. The fact that $\varphi^B$ is convex is a consequence of the general fact that the probability measure of a mean-zero Gaussian process is log-concave, see for instance Lemma 1.1 in [4]. □

*Proof of Lemma 4.* First, let us check that for any $h$ in the RKHS $\mathbb{H}$ of $X_t^\alpha$, it holds $\|h\|_\mathbb{H}^2 = h(0)^2 + \|h - h(0)\|_\mathbb{H}^2$. We use the well-known fact that if a process $X$ is a sum of two independent centered Gaussian components $V$ and $W$, with supports $\mathbb{B}^V$ and $\mathbb{B}^W$ and RKHS $\mathbb{H}^V$ and $\mathbb{H}^W$ respectively, such that $\mathbb{B}^V \cap \mathbb{B}^W = \{0\}$ and $\mathbb{B}^V$ is complemented by a closed subspace that contains $\mathbb{B}^W$, then the RKHS $\mathbb{H}$ of $X$ is the direct sum of $\mathbb{H}^V$ and $\mathbb{H}^W$ and $\|h^V + h^W\|_\mathbb{H}^2 = \|h^V\|_{\mathbb{H}^V}^2 + \|h^W\|_{\mathbb{H}^W}^2$, see for instance Lemma 9.1 in [15]. We apply this fact to the decomposition $X_t^\alpha = V + W$, with $V = Z_0$ and $W = X_t^\alpha - Z_0$, see Equation (9). The support $\mathbb{B}^V$ is the set of all constant functions, while $\mathbb{B}^W$ is included in the (closed) set of all continuous functions $f$ with $f(0) = 0$. Since $\mathbb{B}^V \cap \mathbb{B}^W = \{0\}$, the preceding result implies the announced decomposition.

Now note that

$$\inf_{h \in \mathbb{H},\ \|h - w_0 - \rho\|_\infty < \varepsilon} \|h\|_\mathbb{H}^2 = \inf_{g \in \mathbb{H},\ \|g - w_0\|_\infty < \varepsilon} \|g + \rho\|_\mathbb{H}^2.$$

For any $g$ belonging to the set defining the latter infimum,

$$\begin{aligned} \|g+\rho\|_{\mathbb{H}}^2 &= \|g+\rho-g(0)-\rho\|_{\mathbb{H}}^2 + (g(0)+\rho)^2 \\ &= \|g-g(0)\|_{\mathbb{H}}^2 + (g(0)+\rho)^2 = \|g\|_{\mathbb{H}}^2 + 2g(0)\rho + \rho^2. \end{aligned}$$

Since $\|g-w_0\|_\infty < \varepsilon$ in particular we have $|g(0)-w_0(0)| < \varepsilon$, which gives the desired bound on the infimum and hence on the concentration function. □

*Proof of Lemma 7.* From the proof of Theorem 14 in [7][p. 588], we know that for any $0 \le t \le 1$ and $0 < u \le t$, it holds

$$|I_{0+}^{1-\delta} f(t-u) - I_{0+}^{1-\delta} f(t)| \lesssim u + u \int_1^{t/u} w^\delta \{(w-1)^{-\delta} - w^{-\delta}\} dw.$$

Since $\delta \in (0,1)$, the latter integral is bounded if $t/u \le 2$. If $t/u > 2$ we split the integral in a part over $[1,2]$, which is bounded, and a part over $[2, t/u]$. For the latter part, the mean value theorem gives $|(w-1)^{-\delta} - w^{-\delta}| \le (w-1)^{-\delta-1}$. Thus using the inequality $w \le 2(w-1)$ for $w \ge 2$, we obtain that the integrand is bounded from above by $(w-1)^{-1}$, which leads to

$$|I_{0+}^{1-\delta} f(t-u) - I_{0+}^{1-\delta} f(t)| \lesssim u(1 + \log(1+t/u)) \tag{10}$$

But this also holds for $t \le u$ since then by definition $I_{0+}^{1-\delta} f(t-u) = 0$ and we can use the preceding display with $t=u$ to get that $|I_{0+}^{1-\delta} f(t)| \lesssim t \lesssim u\{1 + \log(1+t/u)\}$. Thus using Fubini's theorem and then (10), one obtains that for any $t > 0$ and any real $u$,

$$\begin{aligned} |I_{0+}^{1-\delta}(f*g)(t)| &\lesssim \int |I_{0+}^{1-\delta} f(t-u) - I_{0+}^{1-\delta} f(t)| |g(u)| du \\ &\lesssim \int |u|\{1 + \log(1+t/|u|)\} |g(u)| du. \end{aligned}$$

Hence by the Cauchy-Schwarz inequality

$$\begin{aligned} \|I_{0+}^{1-\delta}(f*g)\|_2^2 &\lesssim \int_0^1 \left( \int \{1+\log(1+t/|u|)\}^2 u^2 g(u)^2 du \right) dt \\ &\lesssim \int u^2 \{1 + \log^2(1+|u|^{-1})\} g(u)^2 du. \quad □ \end{aligned}$$

## 5. Conclusion

We have defined a notion of lower bound for the rate of convergence of the posterior distribution and given a scheme to obtain lower bounds in a nonparametric framework when the prior is a Gaussian process. Lower and upper bound rates turn out to be intimately related to the behavior of the concentration function $\varphi_{f_0}$ of the Gaussian process at the true $f_0$. When $f_0$ is smooth enough, the small ball term in $\varphi_{f_0}$ dominates and determines the rate. On the contrary, when the prior is much smoother than the function, the RKHS-approximation term dominates and in general some extra information on $f_0$ is needed in order to determine the precise behavior of $\varphi_{f_0}$ explicitly. In the framework of Section 3.1 we were able to obtain that known upper bound rates are, up to constants or log factors, also lower bounds rates, thus leading to optimality of these rates up to constants or log factors. In Section 3.2 we have obtained lower bound results for the posterior rate when the prior is itself non-Gaussian (though constructed from a Gaussian prior) using Lemma 1 directly. Since the proof of Theorem 1 on Gaussian priors also relies on this result, our work also underlines the usefulness of Lemma 1 in obtaining lower bound results.

## Acknowledgments

The author would like to thank the Associate Editor and the referees for their detailed and useful comments.